# On the Duality Theory for Problems with Higher Order Differential Inclusions

Elimhan N. Mahmudov[a,b]

[a]Department of Mathematics, Istanbul Technical University, 34469 Maslak, Istanbul, Turkey,
[b]Azerbaijan National Academy of Sciences Institute of Control Systems, Azerbaijan. e-mail: elimhan22@yahoo.com; Tel: +90 (212) 285 3253, orcid.org/0000-0003-2879-6154

**Abstract**. This paper on the whole concerns with the duality of Mayer problem for $\kappa$-th order differential inclusions, where $\kappa$ is an arbitrary natural number. Thus, this work for constructing the dual problems to differential inclusions of any order can make a great contribution to the modern development of optimal control theory. To this end in the form of Euler-Lagrange type inclusions and transversality conditions the sufficient optimality conditions are derived. The principal idea of obtaining optimal conditions is locally adjoint mappings. It appears that the Euler-Lagrange type inclusions for both prımary and dual problems are "duality relations". To demonstrate this approach, some semilinear problems with $\kappa$-th order differential inclusions are considered. Also, the optimality conditions and the duality theorem in problems with second order polyhedral differential inclusions are proved. These problems show that maximization in the dual problems are realized over the set of solutions of the Euler-Lagrange type differential inclusions/equations.

## 1. Introduction

It is well known that many extremal problems, such as classical problems of optimal control, differential games, models of economic dynamics, macroeconomic problems, etc. [[are described in terms of set-valued mappings and form a component part of the modern mathematical theory of controlled dynamical systems [2,6,8-12,14,17,36,37 . In the papers [19,20, 22-30] and monographs [21,33] for different optimal control problems of discrete processes and differential inclusions (DFIs) are formulated necessary and sufficient optimality conditions. Along with these duality is a powerful and widely employed tool in convex optimization theory for a number of reasons and the convexity properties of the primal problem are important in establishing the connection between it and the dual problem. First, the dual problem is always concave even if the primal is not convex. Second, the number of variables in the dual problem is equal to the number of constraints in the primal which is often less than the number of variables in the primal problem. Third, the maximum value achieved by the dual problem is often equal to the minimum of the primal (see, for example, [7,9,15,18,34,35] and references therein). In [35] is provided a beautiful framework for the analysis of the optimality conditions of rather abstract mathematical programming problems in terms of dual program. An emphasis is placed on the fundamental importance of the concepts of Lagrangian function, saddle-point, and saddle-value. General examples are drawn from nonlinear programming, approximation, stochastic programming, the calculus of variations, and optimal control. The paper [34] on the whole using Fenchel conjugates develops a geometric approach of variational analysis for the case of convex objects considered in locally convex topological spaces and also in Banach space settings.

Besides, in the field of differential inclusions of the higher order, serious qualitative research is intensively carried out [3-5,31,32]. The aim of the papers [4,31] is to establish the existence results of viable solutions for a higher-order differential inclusions. The paper [5] studies the three-point boundary value problems for second-order perturbed differential inclusions. The existence of solutions is proved under the condition of non-convexity for a multivalued mapping.

Along with the Lagrangian and Fenchel dualities [1,7,13,15,16], in some papers the concept of infimal convolution of convex functions are successfully used [7,13,18,20,35]. In the work [7] are given two generalized Moreau–Rockafellar-type results for the sum of a convex function with a composition of convex functions in separated locally convex spaces. Then the stable strong duality for composite convex optimization problems is equivalently

characterized

In the present work, the optimality conditions for a $\kappa$-th order DFIs together with their duality approach were considered for the first time. To the best of our knowledge, there are a few papers [34,35] devoted to duality problems of first order DFIs. We investigate the dual results in accordance with the dual operations of addition and infimal convolution of convex functions. Apparently, the difficulties encountered are due to the fact that this approach requires the construction of the duality of coupled discrete and discrete-approximate problems.

Thus, the present paper is dedicated to one of the difficult and interesting fields − duality of optimization problems with any $\kappa$-th order ordinary DFIs. The posed problem and its duality are new. The paper is organized in the following order:

In *Section 2*, the needed facts and supplementary results from the book of Mahmudov [21] are given; Hamiltonian function $H$ and argmaximum sets of a set-valued mapping $F$, the locally adjoint mapping (LAM), conjugate function for Hamiltonian function taken with a minus sign are introduced and the problems for higher order DFIs ($P_H$) with initial point constraints are formulated.

In *Section 3* sufficient conditions of optimality for $\kappa$-th order DFIs are proven. An adjoint inclusion, which naturally is a generalization of the Euler-Lagrange-type inclusion, and associated transversality conditions at the endpoints $t=0$, $t=1$ are formulated. Then a problem ($P_H$) with semilinear higher order DFIs is considered.

In the paper [25], is considered the similar Bolza problem with initial value problem of Cauchy; the main goal here is to obtain sufficient optimality conditions for the Cauchy problem for higher-order differential inclusions. The sufficient conditions including distinctive transversality condition are proved incorporating the Euler–Lagrange and Hamiltonian type inclusions.

*Section 4* deals with the duality problem $(P_H^*)$ of the $\kappa$-th order DFIs, which contains the main results of this paper. Here we treat dual results according to the dual operations of addition and infimal convolution of convex functions. But the construction of the duality problem would lead us too far astray from the main themes of this paper and is therefore omitted. And in this sense the obtained results of Section 4 are only the visible part of the "icebergs". Further, it is shown that a family of functions $\{\tilde{x}^*(t), \tilde{\eta}_k^*(t), k=1,...,\kappa-1\}$ is an optimal solution of the dual problem $(P_H^*)$ if and only if the conditions (*i*)–(*iii*) of Theorem 3.1 are satisfied. In addition, the optimal values in the primary convex $(P_H)$ and the dual concave $(P_H^*)$ problems are equal: $\inf(P_H) = \sup(P_H^*)$. Thus, it is proved that the Euler-Lagrange type adjoint inclusion at the same time is a dual relation, that is a pair of solutions of primary and dual problems satisfies this relation. At the end of this section is considered a Mayer problem with so-called semilinear $\kappa$-th order DFIs. It is obvious that this method, which is certainly of independent interest from qualitative viewpoint, can play an important role in numerical procedures as well. We believe that relying to the method described in this paper it can be obtained the similar duality results to optimal control problems with any higher order differential inclusions.

*In Section 5* optimality conditions and duality theorem in problems with second order polyhedral DFIs are proved.

In this work we pursue a twofold goal. First, we have proven the duality theorem for a problem to continuous problem $(P_H)$. Second, we use this method to establish a duality relation problem to a continuous Mayer problem.

**2. *Necessary facts and problem statement***

The needed concepts, definitions and notions can be found in the book of Mahmudov [21]. Let $\mathbb{R}^n$ be a $n$-dimensional Euclidean space, $(x,v)$ be a pair of elements $x,v \in \mathbb{R}^n$, $\langle x,v \rangle$ be an inner product of $x$ and $v$. Assume that $F:(\mathbb{R}^n)^\kappa \rightrightarrows \mathbb{R}^n$ is a set-valued (multivalued) mapping from $(\mathbb{R}^n)^\kappa = \underbrace{\mathbb{R}^n \times \mathbb{R}^n \times \cdots \times \mathbb{R}^n}_{\kappa}$ into the family of subsets of $\mathbb{R}^n$. Then $F$ is convex if its $\operatorname{gph} F = \{(x,v_1,...,v_{\kappa-1},v): v \in F(x,v_1,...,v_{\kappa-1})\}$ is a convex subset of $(\mathbb{R}^n)^{\kappa+1}$. The mapping $F$ is convex closed if its graph is a convex closed set in $(\mathbb{R}^n)^{\kappa+1}$. $F$ is convex-valued if $F(x,v_1,...,v_{\kappa-1})$ is a convex set for each $(x,v_1,...,v_{\kappa-1}) \in \operatorname{dom} F = \{(x,v_1,...,v_{\kappa-1}) : F(x,v_1,...,v_{\kappa-1}) \neq \varnothing\}$. The Hamiltonian function and argmaximum set for a set-valued mapping $F$ are defined as

$$H_F(x,v_1,...,v_{\kappa-1},v^*) = \sup_v \{\langle v,v^* \rangle : v \in F(x,v_1,...,v_{\kappa-1})\}, \ v^* \in \mathbb{R}^n,$$

$$F_A(x,v_1,...,v_{\kappa-1};v^*) = \{v \in F(x,v_1,...,v_{\kappa-1}) : \langle v,v^* \rangle = H_F(x,v_1,...,v_{\kappa-1},v^*)\},$$

respectively. For a convex $F$ we put $H_F(x,v_1,...,v_{\kappa-1},v^*) = -\infty$ if $F(x,v_1,...,v_{\kappa-1}) = \varnothing$.

A convex cone $K_A(z)$, $z=(x,v_1,...,v_{\kappa-1},v)$ is a cone of tangent directions at a point $z \in A \subset (\mathbb{R}^n)^{\kappa+1}$ if from $\bar{z}=(\bar{x},\bar{v}_1,...,\bar{v}_{\kappa-1},\bar{v}) \in K_A(z)$ it follows that $\bar{z}$ is a tangent vector at point $z \in A$, in other words if there exists a function $\varphi(\lambda) \in (\mathbb{R}^n)^{\kappa+1}$ such that $z + \lambda \bar{z} + \varphi(\lambda) \in A$ for sufficiently small $\lambda > 0$ and $\lambda^{-1}\varphi(\lambda) \to 0$, as $\lambda \downarrow 0$.

For a such mapping $F$ at a point $(x^0, v_1^0, ..., v_{\kappa-1}^0, v^0) \in \operatorname{gph} F$ a cone of tangent directions is defined as follows

$$K_{\mathrm{gph}F}(x^0, v_1^0, ..., v_{\kappa-1}^0, v^0) = \mathrm{cone}\left[\mathrm{gph}\, F - (x^0, v_1^0, ..., v_{\kappa-1}^0, v^0)\right]$$

$$= \left\{(\bar{x}, \bar{v}_1, ..., \bar{v}_{\kappa-1}, \bar{v}) : \bar{x} = \lambda(x - x^0),\ \bar{v}_j = \lambda(v_j - v_j^0),\ \bar{v} = \lambda(v - v^0),\ j = 1,.., \kappa-1\right\},$$

$$\forall (x, v_1, v_2, ... v_\kappa) \in \mathrm{gph}F.$$

For a convex mapping $F$ a set-valued mapping $F^*(\cdot, x, v_1, ..., v_{\kappa-1}, v) : \mathbb{R}^n \rightrightarrows (\mathbb{R}^n)^\kappa$ defined by

$$F^*\left(v_\kappa^*; (x^0, v_1^0, ..., v_{\kappa-1}^0, v^0)\right) = \left\{(x^*, v_1^*, ..., v_{\kappa-1}^*) : (x^*, v_1^*, ..., v_{\kappa-1}^*, -v^*) \in K^*_{\mathrm{gph}F}(x^0, v_1^0, ..., v_{\kappa-1}^0, v^0)\right\}$$

is called a locally adjoint mapping (LAM) to $F$ at a point $(x^0, v_1^0, ..., v_{\kappa-1}^0, v^0) \in \mathrm{gph}F$, where $K^*_{\mathrm{gph}F}(x^0, v_1^0, ..., v^0)$ is the dual to a cone of tangent vectors $K_{\mathrm{gph}F}(x^0, v_1^0, ..., v_{\kappa-1}^0, v^0)$.

A set-valued mapping defined by

$$F^*(v^*; (x^0, v_1^0, ..., v_{\kappa-1}^0, v^0)) := \Big\{(x^*, v_1^*, ... v_{\kappa-1}^*) : H_F(x, v_1, v_2, ... v_{\kappa-1}, v^*) - H_F(x^0, v_1^0, ..., v_{\kappa-1}^0, v^*)$$

$$\leq \langle x^*, x - x^0 \rangle + \sum_{j=1}^{\kappa-1} \langle v_j^*, v_j - v_j^0 \rangle, \forall (x, v_1, v_2, ... v_{\kappa-1}) \in (\mathbb{R}^n)^\kappa \Big\},\ v \in F(x, v_1, v_2, ... v_{\kappa-1}; v^*)$$

is called the LAM to "nonconvex" mapping $F$ at a point $(x^0, v_1^0, ..., v_{\kappa-1}^0, v^0) \in \mathrm{gph}F$. Obviously, in the convex case $H_F(\cdot, v_\kappa^*)$ is concave and the latter definition of LAM coincide with the previous definition of LAM.

A function $\varphi = \varphi(x, v_1, ..., v_{\kappa-1})$ is called a proper function if it does not assume the value $-\infty$ and is not identically equal to $+\infty$. Obviously, $\varphi$ is proper if and only if $\mathrm{dom}\varphi \neq \varnothing$ and $\varphi(x, v_1, ..., v_{\kappa-1})$ is finite for $(x, v_1, ..., v_{\kappa-1}) \in \mathrm{dom}\varphi = \{(x, v_1, ..., v_{\kappa-1}) : \varphi(x, v_1, ..., v_{\kappa-1}) < +\infty\}$.

**Definition 2.1** A function $\varphi(x, v_1, ..., v_{\kappa-1})$ is said to be a closure if its epigraph $\operatorname{epi}\varphi = \{(\xi, x, v_1, ..., v_{\kappa-1}): \xi \geq \varphi(x, v_1, ..., v_{\kappa-1})\}$ is a closed set.

**Definition 2.2** The function $\varphi^*(x^*, v_1^*, ..., v_{\kappa-1}^*)$ defined as below is called the conjugate of $\varphi$:

$$\varphi^*(x^*, v_1^*, ..., v_{\kappa-1}^*) = \sup_{x, v_1, ..., v_{\kappa-1}} \{\langle x, x^*\rangle + \langle v_1, v_1^*\rangle + ... + \langle v_{\kappa-1}, v_{\kappa-1}^*\rangle - \varphi(x, v_1, ..., v_{\kappa-1})\}.$$

It is clear to see that the conjugate function is closed and convex.

Let us denote

$$M_G(x^*, v_1^*, ..., v_{\kappa-1}^*, v^*) = \inf_{x, v_1, ..., v_{\kappa-1}, v} \{\langle x, x^*\rangle + \langle v_1, v_1^*\rangle + ...$$
$$+ \langle v_{\kappa-1}, v_{\kappa-1}^*\rangle - \langle v, v^*\rangle : (x, v_1, ..., v_{\kappa-1}, v) \in \operatorname{gph}F\}$$

that is, for every $(x, v_1, ..., v_{\kappa-1}) \in (\mathbb{R}^n)^\kappa$

$$M_F(x^*, v_1^*, ..., v_{\kappa-1}^*, v^*) \leq \langle x, x^*\rangle + \langle v_1, v_1^*\rangle + ... + \langle v_{\kappa-1}, v_{\kappa-1}^*\rangle - H_F(x, v_1, ..., v_{\kappa-1}, v^*).$$

It is easy to see that the function

$$M_F(x^*, v_1^*, ..., v_{\kappa-1}^*, v^*) = \inf_{x, v_1, ..., v_{\kappa-1}} \{\langle x, x^*\rangle + \langle v_1, v_1^*\rangle + ... + \langle v_{\kappa-1}, v_{\kappa-1}^*\rangle - H_G(x, v_1, ..., v_{\kappa-1}, v^*)\}$$

is a support function taken with a minus sign. Besides, it follows that for a fixed $v^*$

$$M_F(x^*, v_1^*, ..., v_{\kappa-1}^*, v^*) = -\left[H_F(\cdot, v^*)\right]^*(x^*, v_1^*, ..., v_{\kappa-1}^*)$$

that is, $M_F$ is the conjugate function for $H_F(\cdot, v^*)$ taken with a minus sign. By Lemma 2.6 [21, p.64] it is noteworthy to see that $(x^*, v_1^*, ..., v_{\kappa-1}^*)$ is an element of the LAM $F^*$ if and only if

$$M_F(x^*, v_1^*, \ldots, v_{\kappa-1}^*, v^*) = \langle x, x^* \rangle + \langle v_1, v_1^* \rangle + \ldots + \langle v_{\kappa-1}, v_{\kappa-1}^* \rangle - H_G(x, v_1, \ldots, v_{\kappa-1}, v^*).$$

In what follows the main problem of this paper is the convex optimal control problem for higher order DFIs with free endpoint:

$$\text{infimum} \quad \varphi\big(x(1), x'(1), \ldots, x^{(\kappa-1)}(1)\big), \tag{1}$$

$$(P_H) \quad \frac{d^\kappa x(t)}{dt^\kappa} \in F\big(x(t), x'(t), \ldots, x^{(\kappa-1)}(t), t\big), \text{ a.e. } t \in [0,1], \tag{2}$$

$$x(0) \in Q_0, \ x'(0) \in Q_1, \ x''(0) \in Q_2, \ldots, x^{(\kappa-1)}(0) \in Q_{\kappa-1}, \tag{3}$$

where $F(\cdot, t): (\mathbb{R}^n)^\kappa \rightrightarrows \mathbb{R}^n$ is convex set-valued mapping, $\varphi: (\mathbb{R}^n)^\kappa \to \mathbb{R}^1$ is proper convex function and $Q_j \subseteq \mathbb{R}^n$, $j = 0, 1, \ldots, \kappa-1$ are convex subsets, $\kappa$ is an arbitrary fixed natural number. It is required to find an arc $\tilde{x}(t)$ of the problem (1) – (3) for the $\kappa$-th order differential inclusions satisfying (2) almost everywhere (a.e.) on a time interval $[0,1]$ and the initial point constraints (3) that minimizes the Mayer type cost functional $\varphi$. We label this problem as $(P_H)$. A feasible trajectory $x(\cdot)$ is absolutely continuous function on $[0,1]$ together with the higher order derivatives until $\kappa - 1$, for which $x^{(\kappa)}(\cdot) \in L_1^n([0,1])$. Obviously, such class of functions is a Banach space, endowed with the different equivalent norms. Note that for the Cauchy problem, where each set consists of a single point, i.e. $Q_j = \{x_j\}, j = 0, \ldots, \kappa-1$, we obtain the problem studied in [25].

## 3. Sufficient Conditions of Optimality for Higher Order DFIs

Although the main problem of the paper is duality results for higher order convex DFIs $(P_H)$, as an auxiliary result, we should formulate sufficient conditions of optimality for this problem. These conditions are more precise since they involve useful forms of the Weierstrass-Pontryagin condition and higher order Euler-Lagrange type adjoint inclusions. In what follows, we show that this condition is a duality relation. In the reviewed results a distinctive feature of the proof of a sufficient condition is carried out in Theorem 3.1.

First of all we associate with the problem $(P_H)$ the following so-called the $\kappa$-th order Euler-Lagrange type differential inclusion

(i) $\left( (-1)^{\kappa} x^{*(\kappa)}(t) + \eta_{\kappa-1}^{*}{}', \ \eta_{\kappa-1}^{*}(t) + \eta_{\kappa-2}^{*}{}'(t), \ldots, \eta_{2}^{*}(t) + \eta_{1}^{*}{}'(t), \ \eta_{1}^{*}(t) \right)$

$\in F^*\left( x^*(t); (\tilde{x}(t), \tilde{x}'(t), \ldots, \tilde{x}^{(\kappa)}(t)), t \right)$, a.e. $t \in [0,1]$

and the transversality conditions at the endpoints $t=0$ and $t=1$, respectively:

(ii) $(-1)^{j+1} x^{*(j)}(0) + \eta_{j}^{*}(0) \in K_{Q_{\kappa-k-1}}^{*}\left( \tilde{x}^{(\kappa-j-1)}(0) \right), \ j=1,\ldots,\kappa-1; \ -x^*(0) \in K_{Q_{\kappa-1}}^{*}\left( \tilde{x}^{(\kappa-1)}(0) \right),$

(iii) $\left( (-1)^{\kappa} x^{*(\kappa-1)}(1) + \eta_{\kappa-1}^{*}(1), \ (-1)^{\kappa-1} x^{*(\kappa-2)}(1) + \eta_{\kappa-2}^{*}(1), \ldots, \ x^{*}{}'(1) + \eta_{1}^{*}(1), \ -x^*(1) \right)$

$\in \partial \varphi\left( \tilde{x}(1), \tilde{x}'(1), \ldots, \tilde{x}^{(\kappa-1)}(1) \right).$

In general, our notation and terminology are consistent with first order differential inclusions (see, for example, Mordukhovich [33], Mahmudov [21]).

Later on we assume that $x^*(t)$, $t \in [0,1]$ is absolutely continuous function with the higher order derivatives until $\kappa - 1$ and $x^{*(\kappa)}(\cdot) \in L_1^n([0,1])$. Moreover, $\eta_j^*(t), \ j=1,\ldots,\kappa-1, \ t \in [0,1]$ are absolutely continuous and $\eta_j^{*'}(\cdot) \in L_1^n([0,1])$, $j=1,\ldots,\kappa-1$.

At last we formulate the condition ensuring that the LAM $F^*$ is nonempty at a given point:

(iv)    $\tilde{x}^{(\kappa)}(t) \in F\left(\tilde{x}(t), \tilde{x}'(t), \ldots, \tilde{x}^{(\kappa-1)}(t); x^*(t), t\right)$, a.e. $t \in [0,1]$.

With this tool we are now ready for the result, which gives sufficient conditions of optimality for higher order DFIs.

**Theorem 3.1.** Let $\varphi : (\mathbb{R}^n)^\kappa \to \mathbb{R}^1$ be continuous and proper convex functions and $F$ be a convex set-valued mapping. Moreover let $Q_j, j = 0,\ldots,\kappa-1$ be convex sets. Then for the optimality of the trajectory $\tilde{x}(\cdot)$ in the problem ($P_H$) it is sufficient that there exist a family of absolutely continuous functions $\{x^*(t), \eta_j^*(t), j = 1,\ldots,\kappa-1\}$, $t \in [0,1]$ satisfying a.e. the $\kappa$-th order Euler-Lagrange differential inclusion (i), (iv) and transversality conditions (ii), (iii) at the endpoints $t = 0$ and $t = 1$.

*Proof.* It is not hard to see that by using the Theorem 2.1 [21] in term of Hamiltonian function from the condition (i) we obtain the following important relation

$$H_F\left(x(t), x'(t), \ldots, x^{(\kappa-1)}(t), x^*(t), t\right) - H_F\left(\tilde{x}(t), \tilde{x}'(t), \ldots, \tilde{x}^{(\kappa-1)}(t), x^*(t), t\right)$$

$$\leq \left\langle (-1)^\kappa x^{*(\kappa)}(t) + \frac{d\eta_{\kappa-1}^*(t)}{dt},\ x(t) - \tilde{x}(t) \right\rangle + \left\langle \eta_{\kappa-1}^*(t) + \frac{d\eta_{\kappa-2}^*(t)}{dt},\ x'(t) - \tilde{x}'(t) \right\rangle \quad (4)$$

$$\left\langle \eta_{\kappa-2}^*(t) + \frac{d\eta_{\kappa-3}^*(t)}{dt},\ x''(t) - \tilde{x}''(t) \right\rangle + \cdots + \left\langle \eta_3^*(t) + \frac{d\eta_2^*(t)}{dt},\ x^{(\kappa-3)}(t) - \tilde{x}^{(\kappa-3)}(t) \right\rangle$$

$$+ \left\langle \eta_2^*(t) + \frac{d\eta_1^*(t)}{dt},\ x^{(\kappa-2)}(t) - \tilde{x}^{(\kappa-2)}(t) \right\rangle + \left\langle \eta_1^*(t),\ x^{(\kappa-1)}(t) - \tilde{x}^{(\kappa-1)}(t) \right\rangle.$$

In turn by using the definition of the Hamiltonian function, (4) can be converted to the inequality

$$\left\langle (-1)^\kappa x^{*(\kappa)}(t),\ x(t) - \tilde{x}(t) \right\rangle - \left\langle x^{(\kappa)}(t) - \tilde{x}^{(\kappa)}(t),\ x^*(t) \right\rangle$$

$$+\sum_{j=1}^{\kappa-1}\frac{d}{dt}\left\langle \eta_j^*(t),\ x^{(\kappa-j-1)}(t)-\tilde{x}^{(\kappa-j-1)}(t)\right\rangle \geq 0. \qquad (5)$$

It is not hard to see that the first term in (5) can be transformed to the following useful relation

$$\left\langle (-1)^\kappa x^{*(\kappa)}(t), x(t)-\tilde{x}(t)\right\rangle - \left\langle x^{(\kappa)}(t)-\tilde{x}^{(\kappa)}(t), x^*(t)\right\rangle = \frac{d}{dt}\left\langle (-1)^\kappa x^{*(\kappa-1)}(t), x(t)-\tilde{x}(t)\right\rangle$$

$$-\frac{d}{dt}\left\langle (-1)^{\kappa-1} x^{*(\kappa-2)}(t), x'(t)-\tilde{x}'(t)\right\rangle + \frac{d}{dt}\left\langle (-1)^{\kappa-2} x^{*(\kappa-3)}(t), x''(t)-\tilde{x}''(t)\right\rangle$$

$$+\cdots - \frac{d}{dt}\left\langle x^*(t), x^{(\kappa-1)}(t)-\tilde{x}^{(\kappa-1)}(t)\right\rangle. \qquad (6)$$

Then taking into account (6) if we integrate the inequality (5) we derive that

$$\int_0^1 \left[\left\langle (-1)^\kappa x^{*(\kappa)}(t), x(t)-\tilde{x}(t)\right\rangle - \left\langle x^{(\kappa)}(t)-\tilde{x}^{(\kappa)}(t), x^*(t)\right\rangle\right] dt$$

$$+\sum_{j=1}^{\kappa-1}\int_0^1 d\left\langle \eta_j^*(t),\ x^{(\kappa-j-1)}(t)-\tilde{x}^{(\kappa-j-1)}(t)\right\rangle = \left\langle x(1)-\tilde{x}(1), (-1)^\kappa x^{*(\kappa-1)}(1)\right\rangle$$

$$+\left\langle x'(1)-\tilde{x}'(1), (-1)^{\kappa-1} x^{*(\kappa-2)}(1)\right\rangle + \left\langle x''(1)-\tilde{x}''(1), (-1)^{\kappa-2} x^{*(\kappa-1)}(1)\right\rangle$$

$$+\cdots - \left\langle x^{(\kappa-1)}(1)-\tilde{x}^{(\kappa-1)}(1), x^*(1)\right\rangle - \left\langle x(0)-\tilde{x}(0), (-1)^\kappa x^{*(\kappa-1)}(0)\right\rangle$$

$$-\left\langle x'(0)-\tilde{x}'(0), (-1)^{\kappa-1} x^{*(\kappa-2)}(0)\right\rangle - \left\langle x''(0)-\tilde{x}''(0), (-1)^{\kappa-2} x^{*(\kappa-1)}(0)\right\rangle$$

$$-\cdots + \left\langle x^{(\kappa-1)}(0)-\tilde{x}^{(\kappa-1)}(0), x^*(0)\right\rangle + \sum_{j=1}^{\kappa-1}\left\langle \eta_j^*(1),\ x^{(\kappa-j-1)}(1)-\tilde{x}^{(\kappa-j-1)}(1)\right\rangle$$

$$-\sum_{j=1}^{\kappa-1}\left\langle \eta_j^*(0),\ x^{(\kappa-j-1)}(0)-\tilde{x}^{(\kappa-j-1)}(0)\right\rangle = \left\langle x^*(0),\ x^{(\kappa-1)}(0)-\tilde{x}^{(\kappa-1)}(0)\right\rangle$$

$$-\sum_{j=1}^{\kappa-1}\left\langle(-1)^{j+1}x^{*(j)}(0)+\eta_j^*(0),\ x^{(\kappa-j-1)}(0)-\tilde{x}^{(\kappa-j-1)}(0)\right\rangle \quad (7)$$

$$+\sum_{j=1}^{\kappa-1}\left\langle(-1)^{j+1}x^{*(j)}(1)+\eta_j^*(1),\ x^{(\kappa-j-1)}(1)-\tilde{x}^{(\kappa-j-1)}(1)\right\rangle-\left\langle x^*(1),\ x^{(\kappa-1)}(1)-\tilde{x}^{(\kappa-1)}(1)\right\rangle\geq 0.$$

On the other hand, by definition of the dual cone the transversality condition (*ii*) at the point $t=0$ for all feasible trajectories we have

$$\sum_{j=1}^{\kappa-1}\left\langle(-1)^{j+1}x^{*(j)}(0)+\eta_j^*(0),\ x^{(\kappa-j-1)}(0)-\tilde{x}^{(\kappa-j-1)}(0)\right\rangle$$

$$-\left\langle x^*(0),\ x^{(\kappa-1)}(0)-\tilde{x}^{(\kappa-1)}(0)\right\rangle\geq 0,\ \forall x^{(\kappa-j-1)}(0)\in Q_{\kappa-k-1},\ j=0,...,\kappa-1$$

and so we deduce from (7) that

$$\sum_{j=1}^{\kappa-1}\left\langle(-1)^{j+1}x^{*(j)}(1)+\eta_j^*(1),\ x^{(\kappa-j-1)}(1)-\tilde{x}^{(\kappa-j-1)}(1)\right\rangle-\left\langle x^*(1),\ x^{(\kappa-1)}(1)-\tilde{x}^{(\kappa-1)}(1)\right\rangle\geq 0. \quad (8)$$

In turn, by the transversality condition (*iii*) at the endpoint $t=1$ for all feasible trajectories we have

$$\varphi\big(x(1),x'(1),\ldots,x^{(\kappa-1)}(1)\big)-\varphi\big(\tilde{x}(1),\tilde{x}'(1),\ldots,\tilde{x}^{(\kappa-1)}(1)\big)$$

$$\geq\sum_{j=1}^{\kappa-1}\left\langle(-1)^{j+1}x^{*(j)}(1)+\eta_j^*(1),\ x^{(\kappa-j-1)}(1)-\tilde{x}^{(\kappa-j-1)}(1)\right\rangle-\left\langle x^*(1),\ x^{(\kappa-1)}(1)-\tilde{x}^{(\kappa-1)}(1)\right\rangle. \quad (9)$$

Finally, from the equalities (8) and (9) we deduce the desired result, i.e. for all feasible solutions

$$\varphi\big(x(1),x'(1),\ldots,x^{(\kappa-1)}(1)\big)-\varphi\big(\tilde{x}(1),\tilde{x}'(1),\ldots,\tilde{x}^{(\kappa-1)}(1)\big)\geq 0$$

i.e. $\varphi(x(1), x'(1), \ldots, x^{(\kappa-1)}(1)) \geq \varphi(\tilde{x}(1), \tilde{x}'(1), \ldots, \tilde{x}^{(\kappa-1)}(1))$, $\forall x(t), t \in [0,1]$ and $\tilde{x}(t), t \in [0,1]$

is optimal. □

The following corollaries show that the adjoint inclusion (*i*) of Theorem 3.1 can be considered as a natural generalization of the Euler-Lagrange inclusion.

**Corollary 3.1** For a Mayer problem with first order differential inclusions ($\kappa = 1$)

$$\text{infimum } \varphi(x(1)) \text{ subject to } \frac{dx(t)}{dt} \in F(x(t), t), \text{ a.e. } t \in [0,1], x(0) \in Q_0$$

the conditions (*i*)-(*iii*) of Theorem 3.1 consist of the following

(i) $\quad -\frac{dx^*(t)}{dt} \in F^*\left(x^*(t); (\tilde{x}(t), \tilde{x}'(t), t)\right)$, (ii) $-x^*(0) \in K^*_{Q_0}(x(0))$, (iii) $-x^*(1) \in \partial \varphi(\tilde{x}(1))$.

*Proof.* In fact, in this case $\eta_j^* \equiv 0 \, (j = 1, \ldots, \kappa - 1)$. By convention, we consider also that $\eta_0^* \equiv 0$, when $\kappa = 1$. Then the conditions of the corollary are an elementary consequence of conditions (*i*)-(*iii*) of Theorem 3.1. □

**Corollary 3.2** For a Mayer problem with the second order differential inclusions ($\kappa = 2$)

$$\text{infimum } \varphi(x(1), x'(1)),$$

$$\frac{d^2 x(t)}{dt^2} \in F(x(t), x'(t), t), \text{ a.e. } t \in [0,1], x(0) \in Q_0, x'(0) \in Q_1$$

the conditions (*i*)-(*iii*) of Theorem 3.1 consists of the following

(i) $\quad \left(x^{*\prime\prime}(t) + \eta_1^{*\prime}, \eta_1^*\right) \in F^*\left(x^*(t); (\tilde{x}(t), \tilde{x}'(t), \tilde{x}''(t)), t\right)$, a.e. $t \in [0,1]$,

(ii) $x^{*\prime}(0) + \eta_1^*(0) \in K_{Q_0}^*(\tilde{x}(0)); \ -x^*(0) \in K_{Q_1}^*(\tilde{x}'(0))$,

(iii) $\left(x^{*\prime}(1) + \eta_1^*(1), \ -x^*(1)\right) \in \partial\varphi(\tilde{x}(1), \tilde{x}'(1))$.

*Proof.* In fact, in this case, since in the relation $(-1)^{j+1} x^{*(j)}(0) + \eta_j^*(0) \in K_{Q_{\kappa-j-1}}^*\left(\tilde{x}^{(\kappa-j-1)}(0)\right)$ $\kappa = 2$, then $k = 1$. The second order Euler-Lagrange type inclusion and transversality condition at a point $t = 1$ are immediately consequence of the conditions (*i*), (*iii*) of Theorem 3.1. □

Consider now the problem ($P_H$) with semilinear higher order DFIs:

$$\text{infimum } \varphi\left(x(1), x'(1), \ldots, x^{(\kappa-1)}(1)\right),$$

($P_{SL}$) $\qquad \dfrac{d^\kappa x(t)}{dt^\kappa} \in F\left(x(t), x'(t), \ldots, x^{(\kappa-1)}(t), t\right), \text{ a.e. } t \in [0,1],$

$$x(0) \in Q_0, \ x'(0) \in Q_1, \ x''(0) \in Q_2, \ldots, x^{(\kappa-1)}(0) \in Q_{\kappa-1},$$

$$F(x, v_1, \ldots, v_{\kappa-1}) = A_0 x + \sum_{j=1}^{\kappa-1} A_j v_j + BU,$$

where $A_j (j = 0, \ldots, \kappa-1)$ and $B$ are $n \times n$ and $n \times r$ matrices, respectively, $\varphi$ convex proper function as before, $U \subseteq \mathbb{R}^r, Q_j \subseteq \mathbb{R}^n (j = 0, \ldots, \kappa-1)$ are convex sets. Obviously, $F(\cdot, t)$: $\left(\mathbb{R}^n\right)^\kappa \rightrightarrows \mathbb{R}^n$ is convex mapping. The problem is to find a controlling parameter $\tilde{u}(t) \in U$ such that the arc $\tilde{x}(t)$ corresponding to it minimizes $\varphi(x(1), \ldots, x^{(\kappa-1)}(1))$. Let us introduce a set-valued mapping of the form $F(x, v_1, \ldots, v_{\kappa-1}) = A_0 x + \sum_{j=1}^{\kappa-1} A_j v_j + BU$.

In order to formulate the sufficient condition for a problem ($P_{SL}$) it remains to construct the LAM to set-valued mapping $F(\cdot,t):(\mathbb{R}^n)^\kappa \rightrightarrows \mathbb{R}^n$. It is easy to see that the Hamiltonian function is

$$H_F(x,v_1,...,v_{\kappa-1},v^*) = \sup_v \{\langle v,v^* \rangle : v \in F(x,v_1,...,v_{\kappa-1})\}$$

$$= \left\langle A_0 x + \sum_{j=1}^{\kappa-1} A_j v_j, v^* \right\rangle + \sup_v \{\langle v,v^* \rangle : v \in BU\} = \langle x, A_0^* v^* \rangle + \sum_{j=1}^{\kappa-1} \langle v_j, A_j^* v^* \rangle + W_U(B^* v^*). \quad (10)$$

Further, from the definition of argmaximum set and (10) we have

$$F_A(x,v_1,...,v_{\kappa-1};v^*) = \left\{ v \in F(x,v_1,...,v_{\kappa-1}) : \langle v,v^* \rangle = \langle x, A_0^* v^* \rangle + \sum_{j=1}^{\kappa-1} \langle v_j, A_j^* v^* \rangle + W_U(B^* v^*) \right\}. \quad (11)$$

On the other hand, by definition of LAM or equivalently, by Theorem 2.1[21] we can deduce that

$$F^*(v^*;(x,v_1,...,v_{\kappa-1},v)) = \begin{cases} \partial_{(x,v_1,...,v_\kappa)} H_F(x,v_1,...,v_{\kappa-1},v^*), & v \in F_A(x,v_1,...,v_{\kappa-1};v^*), \\ \varnothing, & v \in F_A(x,v_1,...,v_{\kappa-1};v^*), \end{cases} \quad (12)$$

where it is understood that

$$\partial_{(x,v_1,...,v_\kappa)} H_F(x,v_1,...,v_{\kappa-1},v^*) = -\partial_{(x,v_1,...,v_\kappa)} \left[ -H_F(x,v_1,...,v_{\kappa-1},v^*) \right].$$

Therefore, taking into account (10)-(12) we have

$$F^*(v^*;(x,v_1,...,v_{\kappa-1},v)) = \begin{cases} (A_0^* v^*, A_1^* v^*, ..., A_{\kappa-1}^* v^*), & \text{if } \langle Bu, v^* \rangle = W_U(B^* v^*), \\ \varnothing, & \text{if } \langle Bu, v^* \rangle \neq W_U(B^* v^*). \end{cases} \quad (13)$$

Then comparing the Euler-Lagrange type differential inclusion (*i*) and formula of LAM (13) we can write

$$(-1)^{\kappa} A_{\kappa}^{*} \frac{d^{\kappa} x^{*}(t)}{dt^{\kappa}} + \frac{d\eta_{\kappa-1}^{*}(t)}{dt} = A_{0}^{*} x^{*}(t), \quad \eta_{\kappa-1}^{*}(t) + \frac{d\eta_{\kappa-2}^{*}(t)}{dt} = A_{1}^{*} x^{*}(t),$$

$$..., \eta_{2}^{*}(t) + \frac{d\eta_{1}^{*}(t)}{dt} = A_{\kappa-2}^{*} x^{*}(t), \quad \eta_{1}^{*}(t) = A_{\kappa-1}^{*} x^{*}(t). \tag{14}$$

By sequentially differentiation of $\eta_k^*(t)$, $j = 1,...,\kappa-1$ and introducing each $\eta_j^*(t)$ in the previous relation in (14) we have

$$\eta_j^*(t) = A_{\kappa-j}^* x^*(t) - A_{\kappa-j+1}^* x^{*'}(t) + \cdots + (-1)^{k-1} A_{\kappa-1}^* x^{*(k-1)}(t), \quad j = 1,...,\kappa-1$$

and consequently, it follows that

$$(-1)^{\kappa} A_{\kappa}^{*} \frac{d^{\kappa} x^{*}(t)}{dt^{\kappa}} = A_{0}^{*} x^{*}(t) - A_{1}^{*} x^{*'}(t) + \cdots + (-1)^{\kappa-1} A_{\kappa-1}^{*} x^{*(\kappa-1)}(t). \tag{15}$$

Then the latter equation is just the Euler-Lagrange type differential inclusion (equation). By the formula (13) along with this it is obvious that the so-called, Weierstrass-Pontryagin maximum principle should be satisfied

$$\langle B\tilde{u}(t), x^*(t) \rangle = \sup_{u \in U} \langle Bu, x^*(t) \rangle. \tag{16}$$

As a result, we have the following theorem

**Theorem 3.2.** Let $\varphi : (\mathbb{R}^n)^{\kappa} \to \mathbb{R}^1$ be continuous and proper convex functions and $F(x, v_1, ..., v_{\kappa-1}) = A_0 x + \sum_{j=1}^{\kappa-1} A_j v_j + BU$ be a convex set-valued mapping. Moreover let $Q_j, j = 0, ..., \kappa-1$ be convex sets. Then for the optimality of the trajectory $\tilde{x}(t)$ corresponding to the controlling parameter $\tilde{u}(t)$ to the problem $(P_{SL})$ it is sufficient that there exists an absolutely continuous functions $x^*(t)$, $t \in [0,1]$ with the higher order derivatives until $\kappa-1$ satisfying a.e. the $\kappa$-th order Euler-Lagrange differential equation (15), the transversality conditions (*ii*), (*iii*) at the endpoints $t=0, t=1$ and the Weierstrass-Pontryagin maximum principle (16)

## 4. The Dual Problem for Convex DFIs

We call the following problem, labelled $(P_H^*)$, the dual problem to the previous continuous convex problem $(P_H)$:

$$(P_H^*) \quad \sup_{\substack{x^*(\cdot), \eta_k^*(\cdot) \\ j=1,\ldots,\kappa-1}} \Big\{ -\varphi^*\Big((-1)^\kappa x^{*(\kappa-1)}(1) + \eta_{\kappa-1}^*(1), (-1)^{\kappa-1} x^{*(\kappa-2)}(1) + \eta_{\kappa-2}^*(1), \ldots, x^{*\prime}(1) + \eta_1^*(1), -x^*(1)\Big)$$

$$+ \int_0^1 M_{F(\cdot,t)}\Big((-1)^\kappa x^{*(\kappa)}(t) + \eta_{\kappa-1}^{*\prime}(t), \ \eta_{\kappa-1}^*(t) + \eta_{\kappa-2}^{*\prime}(t), \ldots, \eta_2^*(t) + \eta_1^{*\prime}(t),$$

$$\eta_1^*(t), \ x^*(t)\Big) dt - \sum_{j=1}^{\kappa-1} W_{Q_{\kappa-j-1}}\Big((-1)^j x^{*(j)}(0) - \eta_j^*(0)\Big) - W_{Q_{\kappa-1}}\Big(x^*(0)\Big) \Big\}.$$

Here maximization is taken over the set of functions $\{x^*(\cdot), \eta_j^*(\cdot), j=1,\ldots,\kappa-1\}$.

In what follows, we assume that $x^*(t)$, $t \in [0,1]$ is absolutely continuous function with the higher order derivatives until $\kappa-1$ and $x^{*(\kappa)}(\cdot) \in L_1^n([0,1])$. Moreover $\eta_j^*(t)$, $j=1,\ldots,\kappa-1$, $t \in [0,1]$ are absolutely continuous and $\eta_j^{*\prime}(\cdot) \in L_1^n([0,1])$, $j=1,\ldots,\kappa-1$.

In the proof of duality theorem, at the same time we show that the duality relation is the Euler-Lagrange type adjoint inclusion.

We are now in a position to establish our duality relations between $(P_H)$ and $(P_H^*)$.

**Theorem 4.1** Let the conditions of Theorem 3.1 be satisfied and $\tilde{x}(\cdot)$ be an optimal solution of the primary problem $(P_H)$ with convex DFI. Then a family of functions $\{\tilde{x}^*(t), \tilde{\eta}_j^*(t),\ j = 1, \ldots, \kappa-1\},\ t \in [0,1]$ is an optimal solution of the dual problem $(P_H^*)$ if and only if the conditions (*i*)–(*iii*) of Theorem 3.1 are satisfied. In addition, the optimal values in the primary $(P_H)$ and dual $(P_H^*)$ problems are equal.

*Proof.* First of all we prove that for all feasible solutions $x(\cdot)$ and dual variables $\{x^*(\cdot), \eta_j^*(\cdot), j = 1, \ldots, \kappa-1\}$ of the primary $(P_H)$ and dual $(P_H^*)$ problems, respectively the inequality holds:

$$\varphi\left(x(1), x'(1), \ldots, x^{(\kappa-1)}(1)\right)$$

$$\geq -\varphi^*\left((-1)^\kappa x^*(1)^{(\kappa-1)} + \eta_{\kappa-1}^*(1), (-1)^{\kappa-1} x^{*(\kappa-2)}(1) + \eta_{\kappa-2}^*(1), \ldots, x^{*'}(1) + \eta_1^*(1), -x^*(1)\right)$$

$$+ \int_0^1 M_{F(\cdot,t)}\left((-1)^\kappa x^{*(\kappa)}(t) + \eta_{\kappa-1}^{*'}(t),\ \eta_{\kappa-1}^*(t) + \eta_{\kappa-2}^{*'}(t), \ldots, \eta_2^*(t) + \eta_1^{*'}(t),\right. \quad (17)$$

$$\left. \eta_1^*(t), x^*(t)\right)dt - \sum_{j=1}^{\kappa-1} W_{Q_{\kappa-j-1}}\left((-1)^j x^{*(j)}(0) - \eta_j^*(0)\right) - W_{Q_{\kappa-1}}\left(x^*(0)\right)\Bigg\}.$$

Indeed, applying the Young's inequality [21], we can write

$$-\varphi^*\left((-1)^\kappa x^*(1)^{(\kappa-1)} + \eta_{\kappa-1}^*(1), (-1)^{\kappa-1} x^{*(\kappa-2)}(1) + \eta_{\kappa-2}^*(1), \ldots, x^{*'}(1) + \eta_1^*(1), -x^*(1)\right)$$

$$\leq \varphi\left(x(1), x'(1), \ldots, x^{(\kappa-1)}(1)\right) - \left\langle x(1), (-1)^\kappa x^{*(\kappa-1)}(1) + \eta_{\kappa-1}^*(1)\right\rangle$$

$$-\left\langle x'(1), (-1)^{\kappa-1} x^{*(\kappa-2)}(1) + \eta_{\kappa-2}^*(1)\right\rangle - \cdots - \left\langle x^{(\kappa-2)}(1), x^{*'}(1) + \eta_1^*(1)\right\rangle + \left\langle x^{(\kappa-1)}(1), x^*(1)\right\rangle. \quad (18)$$

By using the $M_{F(\cdot,t)}$ and Hamiltonian function, we have

$$M_{F(\cdot,t)}\left((-1)^\kappa x^{*(\kappa)}(t)+\eta_{\kappa-1}^{*}{}'(t),\ \eta_{\kappa-1}^{*}(t)+\eta_{\kappa-2}^{*}{}'(t),\ ...,\ \eta_2^{*}(t)+\eta_1^{*}{}'(t),\ \eta_1^{*}(t),\ x^{*}(t)\right)$$

$$\leq \left\langle x(t),(-1)^\kappa x^{*(\kappa)}(t)+\eta_{\kappa-1}^{*}{}'(t)\right\rangle+\left\langle x'(t),\eta_{\kappa-1}^{*}(t)+\eta_{\kappa-2}^{*}{}'(t)\right\rangle+\cdots+\left\langle x^{(\kappa-1)}(t),\eta_1^{*}(t)\right\rangle$$

$$-\left\langle x^{(\kappa)}(t),x^{*}(t)\right\rangle = \frac{d}{dt}\left\langle x(t),\eta_{\kappa-1}^{*}(t)\right\rangle+\frac{d}{dt}\left\langle x'(t),\eta_{\kappa-2}^{*}(t)\right\rangle+\cdots+\frac{d}{dt}\left\langle x^{(\kappa-2)}(t),\eta_1^{*}(t)\right\rangle \quad (19)$$

$$+\left\langle x(t),(-1)^\kappa x^{*(\kappa)}(t)\right\rangle-\left\langle x^{(\kappa)}(t),x^{*}(t)\right\rangle.$$

Moreover, by using definition of support functions for each $j=0,1,...,\kappa-1$ we can write

$$-W_{Q_{\kappa-j-1}}\left((-1)^j x^{*(j)}(0)-\eta_j^{*}(0)\right)\leq -\left\langle x^{(\kappa-j-1)}(0),(-1)^j x^{*(j)}(0)-\eta_j^{*}(0)\right\rangle,$$

$$-W_{Q_{\kappa-1}}\left(x^{*}(0)\right)\leq -\left\langle x^{(\kappa-1)}(0),x^{*}(0)\right\rangle,\ j=1,...,\kappa-1,$$

whereas

$$-\sum_{j=1}^{\kappa-1}W_{Q_{\kappa-j-1}}\left((-1)^j x^{*(j)}(0)-\eta_j^{*}(0)\right)-W_{Q_{\kappa-1}}\left(x^{*}(0)\right)$$

$$\leq \sum_{j=1}^{\kappa-1}\left\langle x^{(\kappa-j-1)}(0),(-1)^{j+1}x^{*(j)}(0)\right\rangle-\left\langle x^{(\kappa-1)}(0),x^{*}(0)\right\rangle+\sum_{j=1}^{\kappa-1}\left\langle x^{(\kappa-j-1)}(0),\eta_j^{*}(0)\right\rangle. \quad (20)$$

On the other hand, according to calculation of higher order derivatives [22] from (19) we have

$$\int_0^1 M_{F(\cdot,t)}\left((-1)^\kappa x^{*(\kappa)}(t)+\eta_{\kappa-1}^{*}{}'(t),\ \eta_{\kappa-1}^{*}(t)+\eta_{\kappa-2}^{*}{}'(t),\ ...,\ \eta_2^{*}(t)+\eta_1^{*}{}'(t),\ \eta_1^{*}(t),\ x^{*}(t)\right)dt$$

$$\leq \int_0^1 d\left\langle x(t),\eta_{\kappa-1}^{*}(t)\right\rangle+\int_0^1 d\left\langle x'(t),\eta_{\kappa-2}^{*}(t)\right\rangle+\cdots+\int_0^1 d\left\langle x^{(\kappa-2)}(t),\eta_1^{*}(t)\right\rangle$$

$$+\int_0^1 \left[\left\langle x(t),(-1)^\kappa x^{*(\kappa)}(t)\right\rangle - \left\langle x^{(\kappa)}(t), x^*(t)\right\rangle\right] dt = \left\langle x(1), \eta^*_{\kappa-1}(1)\right\rangle + \left\langle x'(1), \eta^*_{\kappa-2}(1)\right\rangle \qquad (21)$$

$$+ \cdots + \left\langle x^{(\kappa-2)}(1), \eta^*_1(1)\right\rangle - \left\langle x(0), \eta^*_{\kappa-1}(0)\right\rangle - \left\langle x'(0), \eta^*_{\kappa-2}(0)\right\rangle + \cdots + \left\langle x^{(\kappa-2)}(0), \eta^*_1(0)\right\rangle$$

$$+ \int_0^1 \left[\left\langle x(t),(-1)^\kappa x^{*(\kappa)}(t)\right\rangle - \left\langle x^{(\kappa)}(t), x^*(t)\right\rangle\right] dt.$$

Summing (18) and (21) we derive

$$-\varphi^*\left((-1)^\kappa x^*(1)^{(\kappa-1)} + \eta^*_{\kappa-1}(1), (-1)^{\kappa-1} x^{*(\kappa-2)}(1) + \eta^*_{\kappa-2}(1), \ldots, x^{*\prime}(1) + \eta^*_1(1), -x^*(1)\right)$$

$$+ \int_0^1 M_{F(\cdot,t)}\left((-1)^\kappa x^{*(\kappa)}(t) + \eta^{*\prime}_{\kappa-1}(t),\ \eta^*_{\kappa-1}(t) + \eta^{*\prime}_{\kappa-2}(t), \ldots,\ \eta^*_2(t) + \eta^{*\prime}_1(t),\ \eta^*_1(t),\ x^*(t)\right) dt$$

$$\leq \varphi\left(x(1), x'(1), \ldots, x^{(\kappa-1)}(1)\right) - \left\langle x(1),(-1)^\kappa x^{*(\kappa-1)}(1)\right\rangle - \left\langle x'(1),(-1)^{\kappa-1} x^{*(\kappa-2)}(1)\right\rangle \qquad (22)$$

$$- \cdots - \left\langle x^{(\kappa-2)}(1), x^{*\prime}(1)\right\rangle + \left\langle x^{(\kappa-1)}(1), x^*(1)\right\rangle - \left\langle x(0), \eta^*_{\kappa-1}(0)\right\rangle - \left\langle x'(0), \eta^*_{\kappa-2}(0)\right\rangle$$

$$+ \cdots + \left\langle x^{(\kappa-2)}(0), \eta^*_1(0)\right\rangle + \int_0^1 \left[\left\langle x(t),(-1)^\kappa x^{*(\kappa)}(t)\right\rangle - \left\langle x^{(\kappa)}(t), x^*(t)\right\rangle\right] dt.$$

Now summing (20) and (22) we deduce that

$$-\varphi^*\left((-1)^\kappa x^*(1)^{(\kappa-1)} + \eta^*_{\kappa-1}(1), (-1)^{\kappa-1} x^{*(\kappa-2)}(1) + \eta^*_{\kappa-2}(1), \ldots, x^{*\prime}(1) + \eta^*_1(1), -x^*(1)\right)$$

$$+ \int_0^1 M_{F(\cdot,t)}\left((-1)^\kappa x^{*(\kappa)}(t) + \eta^{*\prime}_{\kappa-1}(t),\ \eta^*_{\kappa-1}(t) + \eta^{*\prime}_{\kappa-2}(t), \ldots,\ \eta^*_2(t) + \eta^{*\prime}_1(t),\ \eta^*_1(t),\ x^*(t)\right) dt$$

$$-\sum_{j=1}^{\kappa-1} W_{Q_{\kappa-j-1}}\left((-1)^j x^{*(j)}(0) - \eta^*_j(0)\right) - W_{Q_{\kappa-1}}\left(x^*(0)\right) \qquad (23)$$

$$\leq \varphi\left(x(1), x'(1), \ldots, x^{(\kappa-1)}(1)\right) - \left\langle x(1), (-1)^{\kappa} x^{*(\kappa-1)}(1)\right\rangle$$

$$-\left\langle x'(1), (-1)^{\kappa-1} x^{*(\kappa-2)}(1)\right\rangle - \cdots - \left\langle x^{(\kappa-2)}(1), x^{*}{}'(1)\right\rangle + \left\langle x^{(\kappa-1)}(1), x^{*}(1)\right\rangle$$

$$+\sum_{j=0}^{\kappa-1}\left\langle x^{(\kappa-j-1)}(0), (-1)^{j+1} x^{*(j)}(0)\right\rangle + \int_{0}^{1}\left[\left\langle x(t), (-1)^{\kappa} x^{*(\kappa)}(t)\right\rangle - \left\langle x^{(\kappa)}(t), x^{*}(t)\right\rangle\right] dt.$$

To complete the proof of inequality (17) we should proof that the following sum on the right hand side of (23) denoted by $\Omega$ is equal to zero:

$$\Omega = -\left\langle x(1), (-1)^{\kappa} x^{*(\kappa-1)}(1)\right\rangle - \left\langle x'(1), (-1)^{\kappa-1} x^{*(\kappa-2)}(1)\right\rangle - \cdots - \left\langle x^{(\kappa-2)}(1), x^{*}{}'(1)\right\rangle$$

$$+\left\langle x^{(\kappa-2)}(1), x^{*}(1)\right\rangle + \sum_{j=0}^{\kappa-1}\left\langle x^{(\kappa-j-1)}(0), (-1)^{j+1} x^{*(j)}(0)\right\rangle$$

$$+\int_{0}^{1}\left[\left\langle x(t), (-1)^{\kappa} x^{*(\kappa)}(t)\right\rangle - \left\langle x^{(\kappa)}(t), x^{*}(t)\right\rangle\right] dt.$$

To this end let us compute the integral of $\Lambda = \left\langle x(t), (-1)^{\kappa} x^{*(\kappa)}(t)\right\rangle - \left\langle x^{(\kappa)}(t), x^{*}(t)\right\rangle$ separately; it is not hard to see that $\Lambda$ can be transformed as follows

$$\Lambda = \frac{d}{dt}\left\langle x(t), (-1)^{\kappa} x^{*(\kappa-1)}(t)\right\rangle + \frac{d}{dt}\left\langle x'(t), (-1)^{\kappa-1} x^{*(\kappa-2)}(t)\right\rangle$$

$$+\frac{d}{dt}\left\langle x''(t), (-1)^{\kappa-2} x^{*(\kappa-3)}(t)\right\rangle + \cdots + \frac{d}{dt}\left\langle x^{(\kappa-2)}(t), x^{*}{}'(t)\right\rangle - \frac{d}{dt}\left\langle x^{(\kappa-1)}(t), x^{*}(t)\right\rangle. \quad (24)$$

Then in view of (24) we can compute the integral of $\Lambda$ over time interval $[0,1]$ as follows:

$$\int_{0}^{1}\Lambda dt = \left\langle x(1), (-1)^{\kappa} x^{*(\kappa-1)}(1),\right\rangle + \left\langle x'(1), (-1)^{\kappa-1} x^{*(\kappa-2)}(1)\right\rangle$$

$$+\left\langle x''(1),(-1)^{\kappa-2}x^{*(\kappa-3)}(1)\right\rangle+\cdots+\left\langle x^{(\kappa-2)}(1),x^{*\prime}(1)\right\rangle-\left\langle x^{(\kappa-1)}(1),x^{*}(1)\right\rangle$$

$$-\left\langle x(0),(-1)^{\kappa}x^{*(\kappa-1)}(0)\right\rangle-\left\langle x'(0),(-1)^{\kappa-1}x^{*(\kappa-2)}(0)\right\rangle-\left\langle x''(0),(-1)^{\kappa-2}x^{*(\kappa-3)}(0)\right\rangle$$

$$-\cdots-\left\langle x^{(\kappa-2)}(0),x^{*\,\prime}(0)\right\rangle+\left\langle x^{(\kappa-1)}(0),x^{*}(0)\right\rangle.$$

Therefore

$$\Omega=-\left\langle x(1),(-1)^{\kappa}x^{*(\kappa-1)}(1)\right\rangle-\left\langle x'(1),(-1)^{\kappa-1}x^{*(\kappa-2)}(1)\right\rangle-\cdots-\left\langle x^{(\kappa-2)}(1),x^{*\,\prime}(1)\right\rangle$$

$$+\left\langle x^{(\kappa-2)}(1),x^{*}(1)\right\rangle+\sum_{k=0}^{\kappa-1}\left\langle x^{(\kappa-k-1)}(0),(-1)^{k+1}x^{*(k)}(0)\right\rangle+\left\langle x(1),(-1)^{\kappa}x^{*(\kappa-1)}(1)\right\rangle$$

$$+\left\langle x'(1),(-1)^{\kappa-1}x^{*(\kappa-2)}(1)\right\rangle+\left\langle x''(1),(-1)^{\kappa-2}x^{*(\kappa-3)}(1)\right\rangle+\cdots+\left\langle x^{(\kappa-2)}(1),x^{*\,\prime}(1)\right\rangle$$

$$-\left\langle x^{(\kappa-1)}(1),x^{*}(1)\right\rangle-\left\langle x(0),(-1)^{\kappa}x^{*(\kappa-1)}(0)\right\rangle-\left\langle x'(0),(-1)^{\kappa-1}x^{*(\kappa-2)}(0)\right\rangle$$

$$-\left\langle x''(0),(-1)^{\kappa-2}x^{*(\kappa-3)}(0)\right\rangle-\cdots-\left\langle x^{(\kappa-2)}(0),x^{*\,\prime}(0)\right\rangle+\left\langle x^{(\kappa-1)}(0),x^{*}(0)\right\rangle=0.$$

Thus $\Omega=0$ and the inequality (23) or inequality (17) holds.

Furthermore, suppose that a collection of functions $\{\tilde{x}^{*}(\cdot),\tilde{\eta}_{j}^{*}(\cdot),j=1,\ldots,\kappa-1\}$ satisfies the conditions (*i*)–(*iii*) of Theorem 3.1. Then by definition of LAM the Euler-Lagrange type inclusion (*i*) and the condition (*ii*) imply that

$$H_{F}\left(x(t),x'(t),\ldots,x^{(\kappa-1)}(t),\tilde{x}^{*}(t)\right)-H_{F}\left(\tilde{x}(t),\tilde{x}'(t),\ldots,\tilde{x}^{(\kappa-1)}(t),\tilde{x}^{*}(t)\right)$$

$$\leq\left\langle(-1)^{\kappa}\frac{d^{\kappa}\tilde{x}^{*}(t)}{dt^{\kappa}}+\frac{d\tilde{\eta}_{\kappa-1}^{*}(t)}{dt},x(t)-\tilde{x}(t)\right\rangle+\left\langle\tilde{\eta}_{\kappa-1}^{*}(t)+\frac{d\tilde{\eta}_{\kappa-2}^{*}(t)}{dt},x'(t)-\tilde{x}'(t)\right\rangle$$

$$+\cdots+\left\langle \tilde{\eta}_2^*(t)+\frac{d\tilde{\eta}_1^*(t)}{dt}, x^{(\kappa-2)}(t)-\tilde{x}^{(\kappa-2)}(t)\right\rangle+\left\langle \tilde{\eta}_1^*(t), x^{(\kappa-1)}(t)-\tilde{x}^{(\kappa-1)}(t)\right\rangle$$

whence by the definition of function $M_F$ we derive that

$$\left\langle (-1)^\kappa \frac{d^\kappa \tilde{x}^*(t)}{dt^\kappa}+\frac{d\tilde{\eta}_{\kappa-1}^*(t)}{dt}, \tilde{x}(t)\right\rangle+\left\langle \tilde{\eta}_{\kappa-1}^*(t)+\frac{d\tilde{\eta}_{\kappa-2}^*(t)}{dt}, \tilde{x}'(t)\right\rangle$$

$$+\left\langle \tilde{\eta}_2^*(t)+\frac{d\tilde{\eta}_1^*(t)}{dt}, \tilde{x}^{(\kappa-2)}(t)\right\rangle+\left\langle \tilde{\eta}_1^*(t), \tilde{x}^{(\kappa-1)}(t)\right\rangle - H_F\left(\tilde{x}(t),\tilde{x}'(t),\ldots,\tilde{x}^{(\kappa-1)}(t),\tilde{x}^*(t)\right) \quad (25)$$

$$= M_{F(\cdot,t)}\left((-1)^\kappa \tilde{x}^{*(\kappa)}(t)+\tilde{\eta}_{\kappa-1}^{*'}(t),\ \tilde{\eta}_{\kappa-1}^*(t)+\tilde{\eta}_{\kappa-2}^{*'}(t),\ \ldots,\tilde{\eta}_2^*(t)+\tilde{\eta}_1^{*'}(t),\ \tilde{\eta}_1^*(t),\ \tilde{x}^*(t)\right).$$

On the other hand, by the transversality condition (*ii*) at the initial point $t=0$ we have

$$W_{Q_{\kappa-j-1}}\left((-1)^j \tilde{x}^{*(j)}(0)-\tilde{\eta}_j^*(0)\right)=\left\langle \tilde{x}^{(\kappa-j-1)}(0),(-1)^j \tilde{x}^{*(j)}(0)-\tilde{\eta}_j^*(0)\right\rangle, j=0,\ldots,\kappa-1. \quad (26)$$

At last, by Theorem 1.27 [21] the transversality condition (*iii*) is equivalent to the relation

$$\varphi^*\left((-1)^\kappa \tilde{x}^{*(\kappa-1)}(1)+\tilde{\eta}_{\kappa-1}^*(1),(-1)^{\kappa-1}\tilde{x}^{*(\kappa-2)}(1)+\tilde{\eta}_{\kappa-2}^*(1),\ldots,\ \tilde{x}^{*'}(1)+\tilde{\eta}_1^*(1),\ -\tilde{x}^*(1)\right) \quad (27)$$

$$=\left\langle \tilde{x}(1),(-1)^\kappa \tilde{x}^{*(\kappa-1)}(1)+\tilde{\eta}_{\kappa-1}^*(1)\right\rangle+\left\langle \tilde{x}'(1),(-1)^{\kappa-1}\tilde{x}^{*(\kappa-2)}(1)+\tilde{\eta}_{\kappa-2}^*(1)\right\rangle$$

$$+\cdots+\left\langle \tilde{x}^{(\kappa-2)}(1),\tilde{x}^{*'}(1)+\tilde{\eta}_1^*(1)\right\rangle-\left\langle \tilde{x}^{(\kappa-2)}(1),\tilde{x}^*(1)\right\rangle-\varphi\left(\tilde{x}(1),\tilde{x}'(1),\ldots,\tilde{x}^{(\kappa-1)}(1)\right).$$

Then taking into account the relationships (25)-(27) in (23) the inequality sign is replaced by equality and hence for $\tilde{x}(\cdot)$ and $\{\tilde{x}^*(\cdot),\tilde{\eta}_j^*(\cdot), j=1,\ldots,\kappa-1\}$ the equality of values of the primary and dual problems is ensured. $\square$

The following results show that, in particular, for all problems with differential inclusions of arbitrary order from $(P_H^*)$ one can get their dual problems.

**Corollary 4.1** Consider a problem with first order differential inclusions

$$\text{infimum } \varphi(x(1)) \text{ subject to } \frac{dx(t)}{dt} \in F(x(t),t), \text{ a.e. } t \in [0,1], \ x(0) \in Q_0. \qquad (28)$$

Then the dual problem to the primary (28) is the following problem

$$\sup_{x^*(\cdot)} \left\{ -\varphi^*(-x^*(1)) + \int_0^1 M_{F(\cdot,t)}\left(-x^{*\prime}(t), x^*(t)\right) dt - W_{Q_0}(x^*(0)) \right\}.$$

*Proof.* Observing that by convention $\eta_0^*(t) \equiv 0$ ($\kappa = 1$) we get $\varphi^*\left((-1)^\kappa x^{*(\kappa-1)}(1) + \eta_{\kappa-1}^*(1)\right)$

$= \varphi^*(-x^*(1))$ and $M_{F(\cdot,t)}\left((-1)^\kappa x^{*(\kappa)}(t) + \eta_{\kappa-1}^{*\prime}(t), x^*(t)\right) = M_{F(\cdot,t)}\left(-x^{*\prime}(t), x^*(t)\right)$. Besides, since

$\kappa = 1$, then $j = 0$ and $W_{Q_{\kappa-j-1}}\left((-1)^j x^{*(j)}(0) - \eta_j^*(0)\right) = W_{Q_0}(x^*(0))$. Then we have the desired

problem. □

**Corollary 4.2** The dual problem to problem $(P_H)$ with second order differential inclusions

($\kappa = 2$) is the following problem

$$\sup_{x^*(\cdot),\eta_1^*(\cdot)} \left\{ -\varphi^*\left(x^*(1) + \eta_1^*(1), -x^*(1)\right) + \int_0^1 M_{F(\cdot,t)}\left(x^{*\prime\prime}(t) + \eta_1^{*\prime}(t), \eta_1^*(t), x^*(t)\right) dt \right.$$

$$\left. - W_{Q_0}\left(-x^{*\prime}(0) - \psi_1^*(0)\right) - W_{Q_1}(x^*(0)) \right\}.$$

*Proof.* Obviously, in this case, $\eta_j^*(t) \equiv 0, \ j = 2,...,\kappa-1$ and the proof is similar to the one for

Corollary 4.1 and so is omitted. □

Let us construct the dual problem to continuous $(P_{SL})$ with semilinear higher order DFI, where $F(x, v_1, ..., v_{\kappa-1}) = A_0 x + \sum_{j=1}^{\kappa-1} A_j v_j + BU$. Then it is not hard to calculate that

$$M_F(x^*, v_1^*, ..., v_{\kappa-1}^*, v^*) = \inf_{(x, v_1, ..., v_{\kappa-1}, v) \in \text{gph} F} \left\{ \langle x, x^* \rangle + \langle v_1, v_1^* \rangle + ... + \langle v_{\kappa-1}, v_{\kappa-1}^* \rangle - \langle v, v^* \rangle \right\}$$

$$= \inf_{x, v_1, ..., v_{\kappa-1}} \left[ \langle x, x^* - A_0^* v^* \rangle + \langle v_1, v_1^* - A_1^* v^* \rangle + \cdots + \langle v_{\kappa-1}, v_{\kappa-1}^* - A_{\kappa-1}^* v^* \rangle \right] - \sup_{u \in U} \langle u, B^* v^* \rangle$$

$$= \begin{cases} -W_U(B^* v^*), & \text{if } x^* = A_0^* v^*, ..., v_{\kappa-1}^* = A_{\kappa-1}^* v^*, \\ -\infty, & \text{otherwise.} \end{cases}$$

Then according to the dual problem $(P_H^*)$ we can write

$$M_{F(\cdot, t)} \left( (-1)^\kappa x^{*(\kappa)}(t) + \eta_{\kappa-1}^{*'}, \ \eta_{\kappa-1}^*(t) + \eta_{\kappa-2}^{*'}(t), \ ..., \eta_2^*(t) + \eta_1^{*'}(t), \ \eta_1^*(t), \ x^*(t) \right)$$

$$= \begin{cases} -W_U \left( B^* x^*(t) \right), & \text{if } (-1)^\kappa x^{*(\kappa)}(t) + \eta_{\kappa-1}^{*'} = A_0^* x^*(t), ..., \eta_1^*(t) = A_{\kappa-1}^* x^*(t), \\ -\infty, & \text{otherwise.} \end{cases}$$

Hence

$$M_{F(\cdot, t)} \left( (-1)^\kappa x^{*(\kappa)}(t) + \eta_{\kappa-1}^{*'}, \ \eta_{\kappa-1}^*(t) + \eta_{\kappa-2}^{*'}(t), \ ..., \eta_2^*(t) + \eta_1^{*'}(t), \ \eta_1^*(t), \ x^*(t) \right) = -W_U \left( B^* x^*(t) \right),$$

if $(-1)^\kappa \dfrac{d^\kappa x^*(t)}{dt^\kappa} + \dfrac{d\eta_{\kappa-1}^*(t)}{dt} = A_0^* x^*(t), \ \eta_{\kappa-1}^*(t) + \dfrac{d\eta_{\kappa-2}^*(t)}{dt} = A_1^* x^*(t), \ ..., \eta_2^*(t) + \dfrac{d\eta_1^*(t)}{dt} = A_{\kappa-2}^* x^*(t),$

$\eta_1^*(t) = A_{\kappa-1}^* x^*(t)$. By sequentially differentiation of $\eta_j^*(t), \ j = 1, ..., \kappa - 1$ and introducing each $\eta_j^*(t)$ in the previous relation we have

$$\eta_j^*(t) = A_{\kappa-j}^* x^*(t) - A_{\kappa-j+1}^* x^{*'}(t) + \cdots + (-1)^{j-1} A_{\kappa-1}^* x^{*(j-1)}(t), \ j = 1, ..., \kappa - 1 \qquad (29)$$

and consequently, it follows that

$$(-1)^{\kappa}\frac{d^{\kappa}x^{*}(t)}{dt^{\kappa}}=A_{0}^{*}x^{*}(t)-A_{1}^{*}x^{*\prime}(t)+\cdots+(-1)^{\kappa-1}A_{\kappa-1}^{*}x^{*(\kappa-1)}(t). \qquad (30)$$

The latter equation is just the Euler-Lagrange type differential inclusion(equation).

On using (29) for calculating $\eta_j^*(t)$, $j=1,...,\kappa-1$ as is easily verified that the dual problem of problem ($P_{SL}$) is

$$(P_{SL}^*) \quad \sup_{x^*(\cdot)}\left\{-\varphi^*\left((-1)^{\kappa}x^{*(\kappa-1)}(1)+A_1^*x^*(1)-A_2^*x^{*\prime}(1)+\cdots+(-1)^{\kappa-2}A_{\kappa-1}^*x^{*(\kappa-2)}(1),\ (-1)^{\kappa-1}x^{*(\kappa-2)}(1)\right.\right.$$
$$\left.+A_2^*x^*(1)-\cdots+(-1)^{\kappa-1}A_{\kappa-1}^*x^{*(\kappa-3)}(1),...,\ x^{*\prime}(1)+A_{\kappa-1}^*x^*(1),\ -x^*(1)\right)-\int_0^1 W_U(B^*x^*(t))dt$$
$$-\sum_{j=0}^{\kappa-1}W_{Q_{\kappa-j-1}}\left((-1)^j x^{*(j)}(0)-A_{\kappa-j}^*x^*(0)+A_{\kappa-j+1}^*x^{*\prime}(0)-\cdots+(-1)^k A_{\kappa-1}^*x^{*(j-1)}(0)\right)\right\},$$

where $x^*(\cdot)$ is a solution of the adjoint Euler-Lagrange inclusion/equation (30). Consequently, maximization in this dual problem to primary problem ($P_{SL}$) is realized over the set of solutions of the adjoint equation.

**Theorem 4.2** Let the conditions of Theorem 3.1 be satisfied and $\tilde{x}(\cdot)$ be an optimal solution of the primary semilinear problem ($P_{SL}$). Then $\tilde{x}^*(t)$, $t\in[0,1]$ is an optimal solution of the dual problem ($P_{SL}^*$) if and only if the conditions of Theorem 3.2 are satisfied. In addition, the optimal values in the primary ($P_{SL}$) and dual ($P_{SL}^*$) problems are equal.

## 5. Duality in Problems with Second Order Polyhedral DFIs

In this section we establish the dual problem ($PL^*$) to the problem with the following second order polyhedral differential inclusion

$$\text{infimum } \varphi(x(1), x'(1)),$$

$$(PL) \qquad \frac{d^2 x(t)}{dt^2} \in F(x(t), x'(t), t), \text{ a.e. } t \in [0,1],$$

$$x(0) \in Q_0, \ x'(0) \in Q_1, \ F(x, v_1) = \{v : Ax + Bv_1 - Cv \leq d\},$$

where $F$ is a second order polyhedral set-valued mapping, $A, B$ and $C$ are $s \times n$ dimensional matrices, $d$ is a $s$-dimensional column-vector, $\varphi : \mathbb{R}^{2n} \to \mathbb{R}^1$ is a convex function, $Q_0, Q_1$ are nonempty subsets of $\mathbb{R}^n$. We label this problem by (*PL*). Then with respect to the dual problem $(P_H^*)$ first we should calculate $M_F(x^*, v_1^*, v^*)$:

$$M_F(x^*, v_1^*, v^*) = \inf\{\langle x, x^*\rangle + \langle v_1, v_1^*\rangle - \langle v, v^*\rangle : (x, v_1, v) \in \text{gph} F\}. \tag{31}$$

In fact, denoting $w = (x, v_1, v) \in \mathbb{R}^{3n}, w^* = (x^*, v_1^*, -v^*) \in \mathbb{R}^{3n}$ we have a linear programming problem

$$\inf\{\langle w, w^*\rangle : Dw \leq d\}, \ D = [A \vdots B \vdots -C], \tag{32}$$

where $D = [A \vdots B \vdots -C]$ is $s \times 3n$ block matrix. Then according to the linear programming theory if $\tilde{w} = (\tilde{x}, \tilde{v}_1, \tilde{v})$ is a solution of (32), then there exists $s$-dimensional vector $\lambda \geq 0$ such that

$$w^* = -D^*\lambda, \ \langle A\tilde{x} + B\tilde{v}_1 - C\tilde{v} - d, \lambda\rangle = 0.$$

Hence, $w^* = -D^*\lambda$ means that $x^* = -A^*\lambda, \ v_1^* = -B^*\lambda, \ v^* = -C^*\lambda, \ \lambda \geq 0$. Thus, we find that

$$\begin{aligned} M_F(x^*, v_1^*, v^*) &= \langle \tilde{x}, -A^*\lambda\rangle + \langle \tilde{v}_1, -B^*\lambda\rangle - \langle \tilde{v}, -C^*\lambda\rangle \\ &= -\langle A\tilde{x}, \lambda\rangle - \langle B\tilde{v}_1, \lambda\rangle + \langle C\tilde{v}, \lambda\rangle = -\langle d, \lambda\rangle. \end{aligned} \tag{33}$$

On the other hand, from the form of $M_{F(\cdot,t)}(x^{*''}(t)+\eta_1^{*'}(t), \eta_1^*(t), x^*(t))$ of Corollary 4.2 we derive that

$$x^{*''}(t)+\eta_1^{*'}(t)=-A^*\lambda(t),\ \eta_1^*(t)=-B^*\lambda(t),\ x^*(t)=-C^*\lambda(t),\ \lambda(t)\geq 0 \qquad (34)$$

or

$$C^*\lambda''(t)+B^*\lambda'(t)-A^*\lambda(t)=0,\ \lambda(t)\geq 0. \qquad (35)$$

Therefore, taking into account (33)-(35) and Corollary 4.2 we have the following dual problem:

$$\sup_{\lambda(t)\geq 0}\left\{-\varphi^*\left(-C^*\lambda(1)-B^*\lambda(1), C^*\lambda(1)\right)+\int_0^1 M_{F(\cdot,t)}\left(-A^*\lambda(t),-B^*\lambda(t),-C^*\lambda(t)\right)dt\right.$$

$$\left.-W_{Q_0}\left(C^*\lambda'(0)+B^*\lambda(0)\right)-W_{Q_1}\left(-C^*\lambda(0)\right)\right\}.$$

Now, before formulation of duality theorem we should proof sufficient condition of optimality for a problem (*PL*).

**Theorem 5.1.** Let $\varphi:\mathbb{R}^{2n}\to\mathbb{R}^1$ be continuous proper convex function and $F$ be a polyhedral set-valued mapping given in problem (*PL*). Moreover, let $Q_0, Q_1$ be convex sets. Then for the optimality of the trajectory $\tilde{x}(\cdot)$ in problem (*PL*) with second order polyhedral differential inclusions, it is sufficient that there exists a function $\lambda(t)\geq 0,\ t\in[0,1]$ satisfying a.e. the following second order Euler-Lagrange type polyhedral differential inclusion and transversality conditions at the endpoints $t=0$ and $t=1$:

(1) $\quad C^*\lambda''(t)+B^*\lambda'(t)-A^*\lambda(t)=0,\quad \lambda(t)\geq 0$,

$\quad\langle A\tilde{x}(t)+B\tilde{x}'(t)-C\tilde{x}''(t)-d, \lambda(t)\rangle=0,\ \text{a.e.}\ t\in[0,1].$

(2) $\quad -C^*\lambda'(0) - B^*\lambda(0) \in K^*_{Q_0}(\tilde{x}(0));\ C^*\lambda(0) \in K^*_{Q_1}(\tilde{x}'(0)),\ t = 0,$

(3) $\quad \left(-C^*\lambda'(1) - B^*\lambda(1),\ C^*\lambda(1)\right) \in \partial\varphi(\tilde{x}(1), \tilde{x}'(1)),\ t = 1.$

*Proof.* In fact, condition (1) of the theorem follows directly from Theorem 4.1 [26]. Further, by the conditions (ii), (iii) of Corollary 3.2 ($\kappa = 2$) we have

$$x^{*\prime}(0) + \eta_1^*(0) \in K^*_{Q_0}(\tilde{x}(0));\ -x^*(0) \in K^*_{Q_1}(\tilde{x}'(0)), \tag{36}$$

$$\left(x^{*\prime}(1) + \eta_1^*(1),\ -x^*(1)\right) \in \partial\varphi(\tilde{x}(1), \tilde{x}'(1)). \tag{37}$$

Since by (34) $\eta_1^*(t) = -B^*\lambda(t),\ x^*(t) = -C^*\lambda(t)$ we derive from (36), (37) that

$$-C^*\lambda'(0) - B^*\lambda(0) \in K^*_{Q_0}(\tilde{x}(0));\ C^*\lambda(0) \in K^*_{Q_1}(\tilde{x}'(0)),$$

$$\left(-C^*\lambda'(1) - B^*\lambda(1),\ C^*\lambda(1)\right) \in \partial\varphi(\tilde{x}(1), \tilde{x}'(1)).$$

The proof of theorem is completed. $\square$

Then, as a result of Theorems 4.1 and 5.1 for a problem (*PL*) we have the following duality theorem.

**Theorem 5.2** Let the conditions of Theorem 3.1 be satisfied and $\tilde{x}(\cdot)$ be an optimal solution of the primary problem (*PL*). Then $\tilde{x}^*(t), t \in [0,1]$ is an optimal solution of the dual problem (*PL**) if and only if the conditions of Theorem 5.1 are satisfied. In addition, the optimal values in the primary (*PL*) and dual (*PL**) problems are equal.

*6. Conclusion*

The paper deals with the development of Mayer problem for higher order evolution differential inclusions which are often used to describe various processes in science and engineering. First are derived sufficient optimality conditions in the form of Euler-Lagrange type inclusions and transversality conditions. For construction of the duality problem for higher order problem is required skillfully computation of conjugate and support functions. Therefore, to avoid long calculations, construction of duality problem is omitted. It appears that the Euler-Lagrange type inclusions are duality relations for both primary and dual problems. We believe that relying to the unique method described in this paper it can be obtained the similar duality results to optimal control problems with any higher order differential inclusions. The arising difficulties, of course, are connected with the calculation of conjugate function, integral part of dual problem and support functions. There has been a significant development in the study of duality theory to problems with first order differential/difference inclusions in recent years. Besides, there can be no doubt that investigations of duality results to problems with higher order differential inclusions can have great contribution to the modern development of the optimal control theory. Consequently, we can conclude that the proposed method is reliable for solving the various duality problems with higher order differential inclusions.

**Acknowledgment.** In advance, the author wishes to express his sincere thanks to the Editor-in-Chief, Prof. Christiane Tammer for the consideration of this manuscript in the Journal of Optimization.